\documentclass{article}
\usepackage{amsfonts,amssymb}
\usepackage{latexsym}
\begin{document}
\def\dot{\!\cdot\!}
\def\mod{\mathop{\rm mod}\nolimits}
\def\Coker{\mathop{\rm Coker}\nolimits}
\def\Ker{\mathop{\rm Ker}\nolimits}
\def\Re{\mathop{\rm Re}\nolimits}
\def\Im{\mathop{\rm Im}\nolimits}
\def\sh{\mathop{\rm sh}\nolimits}
\def\ch{\mathop{\rm ch}\nolimits}
\def\rk{\mathop{\rm rk}\nolimits}
\def\Cl{\mathop{\rm Cl}\nolimits}
\def\Gal{\mathop{\rm Gal}\nolimits}
\def\kappa{\varkappa}
\def\indlim{\mathop{{\rm ind}\,\lim}}
\def\const{\mathop{\rm const}\nolimits}
\def\Norm{\mathop{\rm Norm}\nolimits}
\def\deg{\mathop{\rm deg}\nolimits}
\def\Hilb{\mathop{\rm Hilb}\nolimits}
\def\Supp{\mathop{\rm Supp}\nolimits}
\def\BS{\mathop{\rm BS}\nolimits}
\def\upper{\mathop{\rm upper}\nolimits}
\def\lower{\mathop{\rm lower}\nolimits}
\def\unc{\mathop{\rm unc}\nolimits}
\def\GL{\mathop{\rm GL}\nolimits}
\def\ab{\mathop{\rm ab}\nolimits}

\font\mathg=eufb7 at 10pt

\title{Asymptotic behaviour of the Euler--Kronecker constant}

\author{{\bf M. A. Tsfasman}}

\date{}

\maketitle

\centerline{\it{To Volodya Drinfeld with friendship and
admiration}}

\vskip 0.3 cm
\footnotetext {Supported in part by the RFBR Grants 02-01-01041,
02-01-22005. }

\noindent {\small{\bf Abstract.}  This appendix to the beautiful
paper \cite{Ih} of Ihara puts it in the context of infinite global
fields of our papers \cite{Ts/Vl} and \cite{Ts/Vl FF}. We study
the behaviour of Euler--Kronecker constant $\gamma_{K}$ when the
discriminant (genus in the function field case) tends to infinity.
Results of \cite{Ts/Vl} easily give us good lower bounds on the
ratio ${{\gamma_{K}}/\log\sqrt{\vert d_{K}\vert}}$. In particular,
for number fields, under the generalized Riemann hypothesis we
prove
$$\liminf{{\gamma_{K}}\over\log\sqrt{\vert d_{K}\vert}}\ge
-0.26049\dots$$
Then we produce examples of class field towers, showing that
$$\liminf{{\gamma_{K}}\over\log\sqrt{\vert d_{K}\vert}}\le
-0.17849\dots$$}

\vskip 0.3 cm


\section{Introduction}

Let $K$ be a global field, i.e., a finite algebraic extension
either of the field ${\mathbb{Q}}$ of rational numbers, or of the
field of rational functions in one variable over a finite field of
constants. Let $\zeta_{K}(s)$ be its zeta-function. Consider its
Laurent expansion at $s=1$
$$\zeta_{K}(s)=c_{-1} (s-1)^{-1}+c_{0}+c_{1}(s-1)+\ldots$$
In \cite{Ih} Yasutaka Ihara introduces and studies the constant
$$\gamma_{K}=c_{0}/c_{-1}.$$ There are several reasons to
study it:

it generalizes the classical Euler constant
$\gamma=\gamma_{\mathbb{Q}} ;$

for imaginary quadratic fields it is expressed by a
beautiful Kronecker limit formula;

for fields with large discriminants its absolute value is at most
of the order of $const \log\sqrt{\vert d_{K}\vert}$, while the
residue $c_{-1}$ itself may happen to be exponential in
$\log\sqrt{\vert d_{K}\vert}$, see \cite{Ts/Vl}.

In this appendix we study asymptotic behaviour of this constant
when the discriminant (genus in the fuction field case) of the
field tends to infinity. It is but natural to compare Ihara's
results \cite{Ih} with the methods of infinite zeta-functions
developed in \cite{Ts/Vl}.

Let $\alpha_{K}=\log\sqrt{\vert d_{K}\vert}$  in the number field
case and $\alpha_{K}=(g_{K}-1)\log q$ in the function field case
over ${\mathbb F}_{q}$. In the number field case Ihara shows that
$$0\ge\limsup_{K}{\gamma_{K}\over\alpha_{K}}\ge
\liminf_{K}{{\gamma_{K}\over\alpha_{K}}\ge -1}.$$ We improve the
lower bound to

\vskip 0.2 cm\noindent {\bf Theorem $1$. }{\em Assuming the
generalized Riemann hypothesis we have \vskip 0.3 cm \noindent
$\liminf\limits_{K}{\gamma_{K}\over\alpha_{K}}\ge$
$$ -{{\log 2 + {1\over 2}\log 3  +{1\over 4}\log 5 +{1\over 6}\log
7}\over{{1\over\sqrt 2 -1}\log 2 + {1\over\sqrt 3 -1}\log 3  +\log
5 +{1\over\sqrt 7 -1}\log 7+{1\over 2}(\gamma+\log{8\pi})}} =
-0.26049\ldots
$$}

\noindent{\bf Remarks.} Unconditionally we get
$\liminf{\gamma_{K}/\alpha_{K}}\ge-0.52227\ldots.$

In the function field case using the same method we get
$0\ge\limsup{\gamma_{K}/\alpha_{K}}\ge
\liminf{\gamma_{K}/\alpha_{K}}\ge -{(\sqrt q +1)^{-1}},$ which, of
course, coincides with Theorem $2$ of Ihara's paper \cite{Ih}.

Let us remark that the upper bound $0$ is attained for any
asymptotically bad family of global fields, and that the lower
bound in the function field case is attained for any
asymptotically {\em optimal} family (such that the ratio of the
number of ${\mathbb F}_{q}$-points to the genus tends to
$\sqrt{q}-1$), which we know to exist whenever $q$ is a square.
Hence, $\limsup{\gamma_{K}/\alpha_{K}}=0$ and in the function
field case with a square $q$, we have $\;\;
\liminf{\gamma_{K}/\alpha_{K}}= -{(\sqrt q +1)^{-1}}$.

\vskip 0.2 cm In Section 3 we construct examples of class field
towers proving (unconditionally)

\vskip 0.2 cm \noindent {\bf Theorem $2$.}
{$$\liminf_{K}\gamma_{K}/\alpha_{K}\le - {{{2\log 2 + {\log 3}
\over{\log \sqrt{5.7.11.13.17.19.23.29.31.37}}}}}
=-0.17849\ldots$$}

This slightly improves the examples given by Ihara in \cite{Ih}.

In the number field case set $\beta_{K}=-({r_{1}\over
2}(\gamma+\log {4\pi})+ {r_{2}}(\gamma+\log {2\pi})).$ If we
complete $\gamma_{K}$ by archimedean terms, we get

\vskip 0.2 cm \noindent {\bf Theorem $3$.} {\em Let
${\tilde{\gamma}_{K}}=\gamma_{K}+\beta_{K}$. Then, under the
generalized Riemann hypothesis, we have
$$\liminf_{K}{{\tilde{\gamma}}_{K}\over\alpha_{K}}\ge
-(\gamma+\log(2\pi))/(\gamma+\log(8\pi))=-0.6353\ldots$$}

It is much easier to see that
$\limsup{{\tilde{\gamma}}_{K}/\alpha_{K}}\le 0$, and that $0$ is
attained for any {\em asymptotically bad} family (i.e., such that
all $\phi$'s vanish, see the definitions below).

The best example we know gives (unconditionally)

\vskip 0.2 cm \noindent {\bf Theorem $4$.}
$$\liminf_{K}{\tilde{\gamma}}_{K}/\alpha_{K}\le -0.5478\ldots$$

\section{Bounds}

Let us consider the asymptotic behaviour of $\gamma_K$. We treat
the number field case (the same argument in the function field
case leads to Theorem 2 of \cite{Ih}). Let $\vert d_K\vert$ tend
to infinity. By Lemma 2.2 of \cite{Ts/Vl} any family of fields
contains an {\em asymptotically exact subfamily}, i.e., such that
for any $q$ there exists the limit $\phi_q$ of the ratio of the
number $\Phi_{q}(K)$ of prime ideals of norm $q$ to the "genus"
$\alpha_K$, and also the limits $\phi_{\mathbb R}$ and
$\phi_{\mathbb C}$ of the ratios of $r_{1}$ and $r_{2}$ to
$\alpha_K$. To find $\liminf{{\gamma}}_{K}/\alpha_{K}$ and
$\liminf{\tilde{\gamma}}_{K}/\alpha_{K}$ it is enough to find
corresponding limits for a given asymptotically exact family, and
then to look for their minimal values. In what follows we consider
only asymptotically exact families.

\vskip 0.2 cm \noindent {\bf Theorem $5$.} {\em For an
asymptotically exact family $\{K\}$ we have
$$\lim_{K} {\gamma_K\over\alpha_K} = -\sum {\phi_q \log q \over
q-1},$$ where $q$ runs over all prime powers.}

\vskip 0.2 cm \noindent{\em Proof.} The right-hand side equals
$\xi^{0}_{\phi}(1)$, where $\xi^{0}_{\phi}(s)$ is the
log-derivative of the infinite zeta-function $\zeta_{\phi}(s)$ of
\cite{Ts/Vl}. The corresponding series converges for $\Re s\ge 1$
(Proposition $4.2$ of \cite{Ts/Vl}). We know (\cite{Ih}, $(1.3.3)$
and $(1.3.4)$) that
$$\gamma_{K}=-\lim_{s\to
1}\left(Z_{K}(s)-{1\over{s-1}}\right),$$ where for $\Re(s)>1$
$$Z_{K}(s)=-{\zeta'_{K}\over\zeta_{K}}(s)=
\sum_{P,k\ge 1}{\log N(P)\over N(P)^{ks}}=
\sum_{q}\Phi_{q}(K){\log q\over {q^{s}-1}}.$$ By the same
Proposition $4.2$,
${\zeta'_{K}\over\zeta_{K}}(s)\to\xi^{0}_{\phi}(s)$ and hence
$\gamma_{K}/\alpha_{K}\to \xi^{0}_{\phi}(1)$.

\vskip 0.2 cm \noindent {\em Proof of Theorem $1$.} We have to
maximize $\sum {\phi_q \log q \over q-1}$ under the conditions:

$\phi_q\ge 0;$

for any prime $p$ we have $\sum\limits_{m=1}^\infty {m
\phi_{p^m}} \le \phi_{\mathbb R} + 2 \phi_{\mathbb C}$;

$ \sum\limits_q { \phi_q \log q \over \sqrt q-1} +
\phi_{\mathbb{R}} (\log2 \sqrt {2\pi}+{\pi \over 4} + {\gamma
\over 2})+\phi_{\mathbb{C}} (\log 8\pi + \gamma)\le 1$ (Basic
Inequality, GRH-Theorem $3.1$ of \cite{Ts/Vl}).

If we put
$$a_{0}=\log\sqrt{8\pi}+{\pi\over 4}+ {\gamma \over 2},
a_{1}=\log 8\pi + \gamma, a_{q}={ \log q \over \sqrt q -1},
b_{0}=b_{1}=0, b_{q}={ \log q \over q-1},$$ we are under
conditions (1)-(4) and (i)-(iv) of Section $8$ of \cite{Ts/Vl}.

Theorem $1$ is now straightforward from Proposition $8.3$ of
\cite{Ts/Vl}. Indeed, the maximum is attained for $\phi_{p^m}=0$
for $m > 1$, $\phi_{\mathbb R}=0,$ and
$\phi_2=\phi_3=\phi_5=\phi_7=2\phi_{\mathbb C}$ (calculation shows
that starting from $p '= 11$ the last inequality of Proposition
8.3 is violated).

\vskip 0.2 cm \noindent {\em Proof of Theorem $3$.} It is much
easier. Since in this case all coefficients are positive and the
ratio of the coefficient of the function we maximize to the
corresponding coefficient of the Basic Inequality is maximal for
$\phi_{\mathbb C}$, the maximum is attained when all $\phi$'s
vanish, except for $\phi_{\mathbb C}$.

\vskip 0.2 cm \noindent {\bf Remarks.} {If we want unconditional
results, then instead of the Basic inequality we have to use
Proposition $3.1$ of \cite{Ts/Vl}:
$$2\sum_{q}\phi_q\log q {\sum_{m=1}^{\infty}{1\over{q^m+1}}}+
\phi_{\mathbb{R}}(\gamma/2+1/2 +\log
2\sqrt\pi)+\phi_{\mathbb{C}}(\gamma+\log 4\pi) \le 1.$$ For
$\tilde\gamma_K/\alpha_K$ one easily gets $$\liminf
{\tilde\gamma_K\over\alpha_K}\ge -{\gamma+\log(2\pi)\over
\gamma+\log(4\pi)}= -0.7770\ldots$$ The calculation for
$\gamma_K/\alpha_K$ is trickier, since the last condition of
Proposition 8.3 of \cite{Ts/Vl} is not violated until very large
primes. Changing the coefficients by the first term
${(q+1)}^{-1}$, Zykin \cite{Zy} gets $$\liminf
{\gamma_K\over\alpha_K}\ge -0.52227\ldots$$ Note that (for an
asymptotically exact family) $1+\tilde\gamma_K/\alpha_K$ is just
the value at $1$ of the log-derivative $\xi(s)$ of the {\em
completed infinite zeta-function} $\tilde\zeta(s)$ of
\cite{Ts/Vl}.}

\section{Examples}

Let us bound $\liminf {\gamma_K/\alpha_K}$ from above. To do this
one should provide some examples of families. The easiest is, just
as in Section $9$ of \cite{Ts/Vl}, to produce quadratic fields
having infinite class field towers with prescribed splitting. The
proof of Theorem $1$ suggests that we should look for towers of
totally complex fields, where 2, 3, 5 and 7 are totally split.
This is however imprecise, because the sum of Proposition $8.3$
varies only slightly when we change $p_{0}$. Therefore, I also
look at the cases when 2, 3, 5, 7 and 11 are split, and when only
2, 3 and 5 are split, or even only 2 and 3. This leads to a slight
improvement on (1.6.30) of \cite{Ih}.

Each of the following fields has an infinite 2-class field tower
with prescribed splitting (just apply Theorem $9.1$ of
\cite{Ts/Vl}), and Theorem $5$ gives the following list.

\noindent For
${\mathbb{Q}}(\sqrt{11.13.17.19.23.29.31.37.41.43.47.53.59.61.67})$
(the example of Thm.9.4 of \cite{Ts/Vl})  $\mathbb R,$ $2$, $3$,
$5$, $7$ totally split, we get $\liminf
{\gamma_K/\alpha_K}\le-0.1515\ldots$



 \noindent For
${\mathbb{Q}}(\sqrt{-13.17.19.23.29.31.37.41.43.47.53.59.61.73.79})$
(the example of Theorem $9.5$ of \cite{Ts/Vl}) with $2$, $3$, $5$,
$7$, and $11$ split we get $-0.1635\ldots$

\noindent For
${\mathbb{Q}}(\sqrt{-7.11.13.17.19.23.29.31.37.41.43.79})$ with
$2$, $3$, $5$ split we get $-0.1727\ldots$

\noindent For
${\mathbb{Q}}(\sqrt{-7.11.13.17.19.23.29.31.37.41.47.59})$ with
$2$, $3$, $5$ split we get $-0.1737\ldots$

\noindent An even better example is found by Zykin \cite{Zy}:
${\mathbb{Q}}(\sqrt{-5.7.11.13.17.19.23.29.31.37})$ with $2$ and
$3$ split gives us $-0.17849\ldots$ This proves Theorem $2$.

For $\liminf {\tilde\gamma_K/\alpha_K}$ the Martinet field
${\mathbb{Q}}(\cos {2\pi\over 11},\sqrt{2},\sqrt{-23})$ (see
Theorem $9.2$ of \cite{Ts/Vl}) gives $-0.5336\ldots$ The best
Hajir-Maire example (see \cite{Ha/Ma}, Section $3.2$) gives
$\liminf {\tilde\gamma_K/\alpha_K}\le -0.5478\ldots$ This proves
Theorem $4$.

{\small{\bf Acknowledgements.} I am grateful to Professor Yasutaka
Ihara for letting me know his results prior to publication and for
a useful e-mail discussion, to my student Alexei Zykin for many
fruitful discussions and for computer verification of my
calculations, and to Don Zagier for several valuable remarks on
the text.}

{\small \noindent Poncelet Laboratory (UMI 2615 of CNRS and the
Independent University of Moscow)

\noindent11 Bolshoi Vlasievskii per., Moscow 119002, Russia; and

\noindent Institut de Math\'ematiques de Luminy, Luminy Case 907,
Marseille 13288, France

\vskip 0.2 cm \noindent E-mail: tsfasman@iml.univ-mrs.fr}


\begin{thebibliography}{99}

\bibitem{Ih} Y. Ihara. On the Euler--Kronecker constants of global
fields and primes with small norms. - This volume.

\bibitem{Ts/Vl} M.A. Tsfasman, S.G. Vl\u{a}du\c{t}. Infinite Global
Fields and the Generalized Brauer--Siegel Theorem. Moscow Math.
J., 2002, v.2, n.2, pp.329-402.

\bibitem{Ts/Vl FF} M.A. Tsfasman, S.G. Vl\u{a}du\c{t}. Asymptotic
properties of zeta-functions. J. Math. Sciences (New York), 1997,
v.84, n.5, pp.1445-1467.

\bibitem{Ha/Ma} F. Hajir, C. Maire.   Tamely ramified towers and
discriminant bounds for number fields II. J. Symbolic Computation,
2002, v.33 , no.4, pp. 415-423.





\bibitem{Zy} A.Zykin. Private communication.

\end{thebibliography}
\end{document}